\documentclass[12pt,fleqn]{article}
\usepackage{mathrsfs}
\usepackage{amsfonts}
\usepackage{amssymb}
\usepackage{amsthm}
\usepackage{enumerate}
\usepackage{latexsym,amsmath}
\usepackage{graphicx}
\usepackage{abstract}
\usepackage{amsmath}
\usepackage{caption}
\title{On the arithmetic-geometric index of graphs}

\date{}

\begin{document}


\author{Shu-Yu Cui$^{a,b}$, Weifan Wang$^b$, Gui-Xian Tian$^b$\footnote{Corresponding author. E-mail address: gxtian@zjnu.cn or guixiantian@163.com.}, Baoyindureng Wu$^c$\\
{\small{\it $^a$Xingzhi College, Zhejiang Normal University,
Jinhua, 321004, China}}\\
{\small{\it $^b$Department of Mathematics, Zhejiang Normal
University, Jinhua, 321004, China}} \\
{\small{\it $^c$College of Mathematics and System Science, Xinjiang University, Urumqi, 830046, China}}
}\maketitle


\begin{abstract}
Very recently, the first geometric-arithmetic index $GA$ and arithmetic-geometric index $AG$ were introduced in mathematical chemistry. In the present paper, we first obtain some lower and upper bounds on $AG$ and characterize the extremal graphs. We also establish various relations between $AG$ and other topological indices, such as the first geometric-arithmetic index $GA$, atom-bond-connectivity index $ABC$, symmetric division deg index $SDD$, chromatic number $\chi$ and so on.
Finally, we present some sufficient conditions of $GA(G)>GA(G-e)$ or $AG(G)>AG(G-e)$ for an edge $e$ of a graph $G$. In particular, for the first geometric-arithmetic index, we also give a refinement of Bollob\'{a}s-Erd\H{o}s-type theorem obtained in [3].



\end{abstract}

\section*{1. Introduction}

\indent \indent We consider only finite, undirected and simple graph throughout this paper. Let $G=(V,E)$ be a simple
graph of order $n$ and size $m$, with vertex set $V(G) = \{{v_1},{v_2}, \cdots {v_n}\} $ and edge set $E(G)$. Denote an edge $e\in
E(G)$ with end vertices $v_i$ and $v_j$ by $v_iv_j$, simply by $i\sim j$. Let ${d_i}$ be the degree of vertex ${v_i}$ for $i=1,2,\ldots,n$. If an edge $e=v_iv_j$ satisfying $d_i=1$, we say that $e$ is a pendent edge and $v_i$ is a pendent vertex. The maximum and minimum degrees of $G$ are denoted by $\Delta $ and $\delta $,
respectively. Let $\overline{d}$ be the average degree of $G$. The minimum non-pendent vertex degree of $G$ is written by $\delta_1$. Also let $p$ denote number of pendent vertices in $G$.

If the vertex set $V(G)$ is the disjoint union of two nonempty subsets $V_1$ and $V_2$, such that every vertex in $V_1$ has
degree $s$ and every vertex in $V_2$ has degree $r$, then $G$ is said to be $(s, r)$\emph{-semiregular}. In particular, if $s =r$, then $G$ is said to be $r$\emph{-regular}. As usual, the complete bipartite graph, the complete graph and the star on $n$ vertices are denoted by ${K_{p,q}}$, ${K_n}$ and ${K_{1,n - 1}}$, respectively.

Topological indices are graph invariant under graph isomorphisms and reflect some
structural properties of the corresponding molecule graph. These indices are found some chemical applications in chemical graph theory, for example, see [4, 8, 9, 10, 13, 17, 18, 19, 20, 21, 22, 29, 30, 32] and the references cited therein. Recently, Vuki\v{c}evi\'{c} and Furtula [27] proposed a newly graph invariant, namely the \emph{first geometric-arithmetic index}, which is defined as follows:
\begin{equation*}
 GA(G) = \sum\limits_{{v_i}{v_j} \in E(G)} {\left( {\frac{{2\sqrt {{d_i}{d_j}} }}{{{d_i} + {d_j}}}} \right)}.
\end{equation*}
They also obtained some bounds on $GA$ index and determined the trees with maximum and minimum $GA$ indices, which are the star and the path, respectively.
In [6], Das  et al. gave some lower and upper bounds on $GA$ index in terms of the order $n$, the size $m$, the minimum degree $\delta$, maximum degree $\Delta$
and the other topological index. In [1], several further inequalities, involving $GA$ index and several other graph parameters, were obtained. Aouchiche and Hansen [2] presented some bounds on $GA$ index in terms of the order $n$, the chromatic number $\chi$, the minimum degree $\delta$, maximum degree $\Delta$ and average degree $\overline{d}$. At the same time, some conjectures were proposed in [2]. Very recently, Chen and Wu [3] disprove four of these conjectures. In addition, they also presented a sufficient condition with $GA(G)>GA(G-e)$ when an edge $e$ is deleted from a graph $G$. For a comprehensive survey and more details on this area, we refer the reader to [7] and references therein.

In 2015, Shegehall and Kanabur [23] introduced the
\emph{arithmetic-geometric index} $AG$ of $G$. It is defined as follows:
\[
AG(G) = \sum\limits_{{v_i}{v_j} \in E(G)} { \frac{1}{2}\left(\sqrt
{\frac{{{d_i}}}{{{d_j}}}}  + \sqrt {\frac{{{d_j}}}{{{d_i}}}}\right)}.
\]
The $AG$ index of path graph with pendant vertices attached to the middle vertices was discussed in [23, 24]. In addition, the
$AG$ index of graphene, which is world's most conductive and effective material for electromagnetic
interference shielding [26], was computed in [25]. Using this newly topological index, Zheng and the present two authors [31] studied spectrum and energy of arithmetic-geometric matrix, in which the sum of all elements is equal to $2AG$. Other bounds of the arithmetic-geometric
energy of graphs were offered in [5,14]. Motivated by these paper, we further consider bounds on the
$AG$ index and discuss the effect on $GA$ and $AG$ indices of deleting an edge from a graph.

For the sake of convenience, here we list other degree-based topological indices, which will be used in subsequent sections.
\begin{itemize}
\item The forgotten index [15] $F\left( G \right) = \sum\limits_{i = 1}^n {d_i^3}
= \sum\limits_{{v_i}{v_j} \in E\left( G \right)} {\left( {d_i^2 + d_j^2} \right)} $.
\item The first Zagreb index [15], ${M_1}(G)= \sum\limits_{i = 1}^n {d_i^2}  = \sum\limits_{_{{v_i}{v_j} \in E(G)}} {\left( {{d_i} + {d_j}} \right)}$.
\item The second Zagreb index [16], ${M_2}(G) = \sum\limits_{_{{v_i}{v_j} \in E(G)}} {{d_i}{d_j}}$.
\item The symmetric division deg index [28], $SDD(G) =\sum\limits_{{v_i}{v_j} \in E(G)} {\left( {\frac{{{d_i}}}{{{d_j}}} + \frac{{{d_j}}}{{{d_i}}}} \right)} $.
\item The atom-bond connectivity index [12], $ABC(G)=\sum\limits_{{v_i}{v_j} \in E(G)} \sqrt{\frac{d_i+d_j-2}{d_id_j}}$.
\end{itemize}

This paper is organized as follows. In Section 2, we present some lower and upper bounds on $AG$ and characterize the extremal graphs. In Section 3, we establish various relations
between $AG$ and other topological indices, such as the first geometric-arithmetic index $GA$, atom-bond-connectivity index $ABC$, symmetric division deg index $SDD$, chromatic number $\chi$ and so on. In Section 4, we obtain some sufficient conditions of $GA(G)>GA(G-e)$ or $AG(G)>AG(G-e)$ for an edge $e$ of a graph $G$. In particular, we give a refinement of Bollob\'{a}s-Erd\H{o}s-type theorem obtained in [3] for the first geometric-arithmetic index. Many examples show that there are considerable differences between $GA$ and $AG$ indices of graphs.

\section*{2. Upper and lower bounds on arithmetic-geometric index}

\indent \indent Throughout this section, we always assume that $G$ is a graph with $p$ pendent vertices. \\

\noindent{\bfseries Theorem 1.} If $G$ is a connected graph of
order $n$ with size $m$, maximum degree $\Delta $, minimum non-pendent vertex degree $\delta_1$, then
\begin{equation}\label{1}
AG(G) \le \frac{{p(\Delta  + 1)}}{{2\sqrt \Delta  }} + \frac{1}{{2\delta _1 }}\sqrt {(m - p)(F + 2M_2  - p(\delta _1  + 1)^2 )} ,
\end{equation}
with equality if and only if $G$ is isomorphic to ${K_{1,n - 1}}$ or $G$ is isomorphic to a regular graph or $G$ is isomorphic to a $(\Delta,1)$-semiregular graph.

\begin{proof}
By the Cauchy-Schwarz inequality, one has
\begin{equation*}
 AG(G) = \sum\limits_{i \sim j} {\frac{1}{2}\left( {\sqrt {\frac{{d_i }}{{d_j }}}  + \sqrt {\frac{{d_j }}{{d_i }}} } \right)}
\end{equation*}
\begin{equation*}
~~~~~~~~~~= \frac{1}{2}\sum\limits_{i \sim j,\;d_j  = 1} {\left( {\sqrt {\frac{{d_i }}{{d_j }}}  + \sqrt {\frac{{d_j }}{{d_i }}} } \right)}  + \frac{1}{2}\sum\limits_{i \sim j,d_i ,d_j  > 1} {\left( {\sqrt {\frac{{d_i }}{{d_j }}}  + \sqrt {\frac{{d_j }}{{d_i }}} } \right)}
\end{equation*}
\begin{equation*}
~~~~~~~~~~\le \frac{{p(\Delta  + 1)}}{{2\sqrt \Delta  }} + \frac{1}{2}\sqrt {(m - p)} \sqrt {\sum\limits_{i \sim j,d_i ,d_j  > 1} {\left( {\sqrt {\frac{{d_i }}{{d_j }}}  + \sqrt {\frac{{d_j }}{{d_i }}} } \right)} ^2 }
\end{equation*}
\begin{equation}\label{2}
~~~~~~~~~~\le \frac{{p(\Delta  + 1)}}{{2\sqrt \Delta  }} + \frac{1}{{2\delta _1 }}\sqrt {(m - p)} \sqrt {\sum\limits_{i \sim j,d_i ,d_j  > 1} {\left( {d_i  + d_j } \right)} ^2 } \end{equation}
\begin{equation*}
~~~~~~~~~~ = \frac{{p(\Delta  + 1)}}{{2\sqrt \Delta  }} + \frac{1}{{2\delta _1 }}\sqrt {(m - p)} \sqrt {\sum\limits_{i \sim j} {\left( {d_i  + d_j } \right)} ^2  - \sum\limits_{i \sim j,d_j  = 1} {\left( {d_i  + d_j } \right)} ^2 }
\end{equation*}
\begin{equation}\label{3}
~~~~~~~~~~ \le \frac{{p(\Delta  + 1)}}{{2\sqrt \Delta  }} + \frac{1}{{2\delta _1 }}\sqrt {(m - p)(F + 2M_2  - p(\delta _1  + 1)^2 )} ~~\text{as}~~ d_i\geq\delta_1,
\end{equation}
implying the required result (\ref{1}).

Now assume that the equality holds in (\ref{1}). Then all inequalities in above proof
must be equalities. From the equality in (\ref{2}), we have ${d_i} =
\Delta$ and $d_j=1$ for any pendent edge $i\sim j$, and $d_i=d_j=\delta_1$ for any non-pendent edge $i\sim j$.
From the equality in (\ref{3}), we have ${d_i} =
\delta_1$ and $d_j=1$ for any pendent edge $i\sim j$. In particular, if $G$ has no pendent edge, that is, $p=0$, then $G$ is isomorphic to a $\Delta$-regular graph. If every edge of $G$ is pendent edge, that is, $m=p$, then $G$ is isomorphic to ${K_{1,n - 1}}$. Otherwise, $0<p<m$, which implies that $G$ is isomorphic to a $(\Delta,1)$-semiregular graph as $G$ is connected.

Conversely, it is easy to check that the equality holds in (\ref{1})
for $K_{1,n-1}$ or a regular graph or a $(\Delta,1)$-semiregular graph.
\end{proof}

\noindent{\bfseries Corollary 1.}
If $G$ is a connected graph of
order $n$ with size $m$, minimum degree $\delta$, then
\begin{equation}\label{1}
AG(G) \le\frac{1}{{2\delta }}\sqrt {m (F + 2M_2 )} ,
\end{equation}
with equality if and only if $G$ is isomorphic to a regular graph.

\begin{proof}
If $G$ has no pendent edge, then $p=0$ and $\delta=\delta_1$. By Theorem 1, we get the required result.
Now assume that $G$ has at least a pendent edge, that is $\delta=1$. We need only to prove that
$
AG(G) \le\frac{1}{{2 }}\sqrt {m (F + 2M_2 )}.
$
Indeed, from the Cauchy-Schwarz inequality, one has
\[
\begin{array}{l}
 (2AG(G))^2  = \left( {\sum\limits_{i \sim j} {\left( {\sqrt {\frac{{d_i }}{{d_j }}}  + \sqrt {\frac{{d_j }}{{d_i }}} } \right)} } \right)^2  \le \left( {\sum\limits_{i \sim j} {(d_i  + d_j )} } \right)^2  \\\\
~~~~~~~~~~~~~~ \le m\sum\limits_{i \sim j} {(d_i  + d_j )^2 }  = m(F + 2M_2 ), \\
 \end{array}
\]
with equality if and only if $d_i=d_j=1$ for a pendent edge $i\sim j$, equivalently, $G$ is isomorphic to $K_2$ as $G$ is connected.
\end{proof}

\noindent{\bfseries Theorem 2.} If $G$ is a connected graph of
order $n$ with size $m$, maximum degree $\Delta $, minimum non-pendent vertex degree $\delta_1$, then
\begin{equation}\label{5}
AG(G) \le \frac{{p(\Delta  + 1)}}{{2\sqrt \Delta  }} + \frac{1}{{2\delta _1 }}\sqrt {F + 2M_2  - p(\delta _1  + 1)^2  + 4\Delta ^2 (m - p)(m - p - 1)} ,
\end{equation}
with equality if and only if $G$ is isomorphic to ${K_{1,n - 1}}$ or $G$ is isomorphic to a regular graph or $G$ is isomorphic to a $(\Delta,1)$-semiregular graph.

\begin{proof}
For the sake of convenience, we first estimate
\[\footnotesize
\begin{array}{l}
 \sum\limits_{i \sim j,d_i ,d_j  > 1} {\left( {\sqrt {\frac{{d_i }}{{d_j }}}  + \sqrt {\frac{{d_j }}{{d_i }}} } \right)}  \\
~~~~ = \sqrt {\sum\limits_{i \sim j,d_i ,d_j  > 1} {\left( {\sqrt {\frac{{d_i }}{{d_j }}}  + \sqrt {\frac{{d_j }}{{d_i }}} } \right)} ^2  + 2\sum\limits_{i \sim j,k \sim l,(v_i,v_j)\neq(v_k ,v_l),d_i ,d_j ,d_k ,d_l  > 1} {\left( {\sqrt {\frac{{d_i }}{{d_j }}}  + \sqrt {\frac{{d_j }}{{d_i }}} } \right)\left( {\sqrt {\frac{{d_k }}{{d_l }}}  + \sqrt {\frac{{d_l }}{{d_k }}} } \right)} }  \\
 \end{array}
\]
\begin{equation}\label{6}\small
~~\le \frac{1}{{\delta _1 }}\sqrt {\sum\limits_{i \sim j,d_i ,d_j  > 1} {\left( {d_i  + d_j } \right)} ^2  + 2\sum\limits_{i \sim j,k \sim l,(v_i,v_j)\neq(v_k ,v_l),d_i ,d_j ,d_k ,d_l  > 1} {\left( {d_i  + d_j } \right)\left( {d_k  + d_l } \right)} }
\end{equation}
\begin{equation}\label{7}
\le \frac{1}{{\delta _1 }}\sqrt {F + 2M_2  - p(\delta _1  + 1)^2  + 8\Delta ^2 \left( {\begin{array}{*{20}c}
   {m - p}  \\
   2  \\
\end{array}} \right)}
\end{equation}
\begin{equation*}
= \frac{1}{{\delta _1 }}\sqrt {F + 2M_2  - p(\delta _1  + 1)^2  + 4\Delta ^2 (m - p)(m - p - 1)}.
\end{equation*}
Note that the function $f(x)=x+\frac{1}{x}$ is an increasing function for $x\geq 1$. Then,
\begin{equation*}
 AG(G)= \frac{1}{2}\sum\limits_{i \sim j,d_j  = 1} {\left( {\sqrt {\frac{{d_i }}{{d_j }}}  + \sqrt {\frac{{d_j }}{{d_i }}} } \right)}  + \frac{1}{2}\sum\limits_{i \sim j,d_i ,d_j  > 1} {\left( {\sqrt {\frac{{d_i }}{{d_j }}}  + \sqrt {\frac{{d_j }}{{d_i }}} } \right)}~~~~~~~~~~~~~~~~~~~~~~~~~~~~~~~~
\end{equation*}
\begin{equation}\label{8}
~~~~\leq\frac{{p(\Delta  + 1)}}{{2\sqrt \Delta  }} + \frac{1}{{2\delta _1 }}\sqrt {F + 2M_2  - p(\delta _1  + 1)^2  + 4\Delta ^2 (m - p)(m - p - 1)}.
\end{equation}

Now assume that the equality holds in (\ref{5}). Then all inequalities in above argument
must be equalities. From the equality in (\ref{6}), we have $d_i=d_j=\delta_1$ for any non-pendent edge $i\sim j$. The equality in (\ref{7}) implies that ${d_i} =
\delta_1$ and $d_j=1$ for any pendent edge $i\sim j$. At the same time, it follows from the equality in (\ref{8}) that ${d_i} =
\Delta$ and $d_j=1$ for any pendent edge $i\sim j$. Observe that $G$ is connected, similar to the argument
of Theorem 1, then $G$ is isomorphic to ${K_{1,n - 1}}$ or $G$ is isomorphic to a regular graph or $G$ is isomorphic to a $(\Delta,1)$-semiregular graph.

Conversely, it is easy to check that the equality holds in (\ref{5})
for $K_{1,n-1}$ or a regular graph or a $(\Delta,1)$-semiregular graph.
\end{proof}

\noindent{\bfseries Corollary 2.}
If $G$ is a connected graph of
order $n$ with size $m$, minimum degree $\delta$ and maximum degree $\Delta $, then
\begin{equation}\label{1}
AG(G) \le\frac{1}{{2\delta }}\sqrt {F + 2M_2+4m(m-1)\Delta^2} ,
\end{equation}
with equality if and only if $G$ is isomorphic to a regular graph.

\begin{proof}
The proof is similar to that of Corollary 1, omitted.
\end{proof}

\noindent{\bfseries Theorem 3.} If $G$ is a connected graph of
order $n$ with size $m$, then
\begin{equation}\label{10}
AG(G) \le \frac{1}{2}p\left( {\frac{n}{{\sqrt {n - 1} }}} \right) + \frac{1}{2}(m - p)\left( {\sqrt {\frac{{n - 1}}{2}}  + \sqrt {\frac{2}{{n - 1}}} } \right),
\end{equation}
with equality if and only if $G$ is isomorphic to a star ${K_{1,n - 1}}$ or $G$ is isomorphic to a complete graph $K_3$.

\begin{proof}
Since the function $f(x)=x+\frac{1}{x}$ is an increasing function for $x\geq 1$. Then, for any pendent edge $i\sim j$ and $d_j=1$,
\begin{equation}\label{11}
\sqrt {\frac{{d_i }}{{d_j }}}  + \sqrt {\frac{{d_j }}{{d_i }}}  \le \sqrt {n - 1}  + \frac{1}{{\sqrt {n - 1} }}=\frac{n}{\sqrt{n-1}},
\end{equation}
with equality if and only if $d_i=n-1$. Similarly, for any non-pendent edge $i\sim j$, one has
\begin{equation}\label{12}
\sqrt {\frac{{d_i }}{{d_j }}}  + \sqrt {\frac{{d_j }}{{d_i }}}  \le \sqrt {\frac{{n - 1}}{2}}  + \sqrt {\frac{2}{{n - 1}}},
\end{equation}
with equality if and only if $d_i=n-1$ and $d_j=2$ for $d_i\geq d_j$.
Therefore,
\[
AG(G) = \frac{1}{2}\sum\limits_{i \sim j,d_j  = 1} {\left( {\sqrt {\frac{{d_i }}{{d_j }}}  + \sqrt {\frac{{d_j }}{{d_i }}} } \right)}  + \frac{1}{2}\sum\limits_{i \sim j,d_i ,d_j  > 1} {\left( {\sqrt {\frac{{d_i }}{{d_j }}}  + \sqrt {\frac{{d_j }}{{d_i }}} } \right)}
\]
\[
~~~~~~~~~ \le \frac{p}{2}\left( {\frac{n}{{\sqrt {n - 1} }}} \right) + \frac{1}{2}(m - p)\left( {\sqrt {\frac{{n - 1}}{2}}  + \sqrt {\frac{2}{{n - 1}}} } \right).
\]

Now assume that the equality holds in (\ref{10}). Then all inequalities in above argument
must be equalities. In the following, without loss of generality, assume that $d_i\geq d_j$ for every edge $i\sim j$.
First if $G$ has no pendent edge, equivalently, $p=0$. Then the equality in (\ref{12}) implies that there is a common neighbor between the end vertices of every edge of $G$. This shows that $G$ is isomorphic to a complete graph $K_3$ as $G$ is connected. Clearly, $G$ is isomorphic to ${K_{1,n - 1}}$ when $p=m$. Finally, assume that $0<p<m$ and ${d_i} =
n-1$, $d_k=1$ for some pendent edge $i\sim k$. Then there must exist a non-pendent edge $i\sim j$ of $G$ such that $d_j=2$ as $m>p$. Thus the vertices $i$ and $j$ have must a common neighbor $l$. Also, from the equality in (\ref{12}), we have $d_l=n-1$. Therefore, $l\sim k$, which implies that $d_k\geq 2$ contradicting to our assumption.

Conversely, it is easy to check that the equality holds in (\ref{10})
for a complete graph $K_3$ or a star $K_{1,n-1}$.
\end{proof}

\noindent{\bfseries Corollary 3.}
If $G$ is a connected graph of
order $n$ with size $m$, minimum degree $\delta\geq2$, then
\begin{equation}\label{13}
AG(G) \le \frac{1}{2}m\left( {\sqrt {\frac{{n - 1}}{2}}  + \sqrt {\frac{2}{{n - 1}}} } \right),
\end{equation}
with equality if and only if $G$ is isomorphic to a complete graph $K_3$.\\

The following lemma comes from [11]. First let $(a_1,a_2,\ldots,a_n)$, $(b_1,b_2,\ldots,b_n)$ be two sequences
of positive real numbers, such that there are positive numbers $A,a,B,b$
satisfying, for any $i\in\{1,2,\ldots,n\}$,
\[
0 < a \le a_i  \le A < \infty ,\;0 < b \le b_i  \le B < \infty.
\]
\noindent{\bfseries Lemma 1}(P\'{o}lya-Szeg\"{o} inequality [11]).
\begin{equation*}
\frac{{\sum\nolimits_{i = 1}^n {a_i^2 } \sum\nolimits_{i = 1}^n
{b_i^2 } }}{{(\sum\nolimits_{i = 1}^n {a_i b_i } )^2 }} \le
\frac{{(ab + AB)^2 }}{{4abAB}},
\end{equation*}
where the equality holds if and only if
\begin{equation*}
p = n \cdot {{\frac{A}{a}} \mathord{\left/
 {\vphantom {{\frac{A}{a}} {(\frac{A}{a} + \frac{B}{b}),q = }}} \right.
 \kern-\nulldelimiterspace} {(\frac{A}{a} + \frac{B}{b}),\; q = }}n \cdot {{\frac{B}{b}} \mathord{\left/
 {\vphantom {{\frac{B}{b}} {(\frac{A}{a} + \frac{B}{b})}}} \right.
 \kern-\nulldelimiterspace} {(\frac{A}{a} + \frac{B}{b})}}
\end{equation*}
are integers and if $p$ of the numbers $a_1,a_2,\ldots,a_n$ are
equal to $a$ and $q$ of these numbers are equal to $A$, and if the
corresponding numbers $b_i$ are equal to $B$ and $b$,
respectively.\\

\noindent{\bfseries Theorem 4.} If $G$ is a connected graph of
order $n$ with size $m$, maximum degree $\Delta $, minimum non-pendent vertex degree $\delta_1$, then
\begin{equation}\label{14}
AG(G)  \ge \frac{{p(\delta _1  + 1)}}{{2\sqrt {\delta _1 } }} + \frac{{\sqrt {2(m - p)(\Delta  + \delta _1 )\sqrt {\Delta \delta _1 } } }}{{\Delta (\Delta  + \delta _1  + 2\sqrt {\Delta \delta _1 } )}}\sqrt {F + 2m\Delta ^2  - p(3\Delta ^2  + 1)},
\end{equation}
with equality if and only if $G$ is isomorphic to ${K_{1,n - 1}}$ or $G$ is isomorphic to a regular graph or $G$ is isomorphic to a $(\Delta,1)$-semiregular graph.

\begin{proof}
For any non-pendent edge $i\sim j$, we get
\begin{equation*}
2 \le \sqrt {\frac{{d_i }}{{d_j }}}  + \sqrt {\frac{{d_j }}{{d_i }}}  \le \sqrt {\frac{\Delta }{{\delta _1 }}}  + \sqrt {\frac{{\delta _1 }}{\Delta }} .
\end{equation*}
Then, take $a=2$, $A=\sqrt {\frac{\Delta }{{\delta _1 }}}  + \sqrt {\frac{{\delta _1 }}{\Delta }}$ and $b=B=1$ in Lemma 1, one has
\begin{equation}\small \label{15}
\left( {\sum\limits_{i \sim j,d_i ,d_j  > 1} {\left( {\sqrt {\frac{{d_i }}{{d_j }}}  + \sqrt {\frac{{d_j }}{{d_i }}} } \right)} } \right)^2  \ge \frac{{8(m - p)(\Delta  + \delta _1 )\sqrt {\Delta \delta _1 } }}{{(\Delta  + \delta _1  + 2\sqrt {\Delta \delta _1 } )^2 }}\sum\limits_{i \sim j,d_i ,d_j  > 1} {\left( {\frac{{d_i }}{{d_j }} + \frac{{d_j }}{{d_i }} + 2} \right)}.
\end{equation}
For the sake of convenience, take
\[
\Gamma_1=\left( {\sum\limits_{i \sim j,d_i ,d_j  > 1} {\left( {\sqrt {\frac{{d_i }}{{d_j }}}  + \sqrt {\frac{{d_j }}{{d_i }}} } \right)} } \right)^2
\]
and
\[
\Gamma_2=\sum\limits_{i \sim j,d_i ,d_j  > 1} {\left( {\frac{{d_i }}{{d_j }} + \frac{{d_j }}{{d_i }} + 2} \right)}.
\]
We first estimate the value of $\Gamma_2$,
\begin{equation*}
\Gamma_2= \sum\limits_{i \sim j,d_i ,d_j  > 1} {\left( {\frac{{d_i }}{{d_j }} + \frac{{d_j }}{{d_i }}} \right)}  + 2(m - p)
\end{equation*}
\begin{equation} \label{16}
~~~\ge \frac{1}{{\Delta ^2 }}\left( {\sum\limits_{i \sim j} {\left( {d_i^2  + d_j^2 } \right)}  - \sum\limits_{i \sim j,d_j  = 1} {\left( {d_i^2  + d_j^2 } \right)} } \right) + 2(m - p)
\end{equation}
\begin{equation}\label{17}
~~~\ge \frac{1}{{\Delta ^2 }}(F - p(\Delta ^2  + 1)) + 2(m - p).
\end{equation}
Plugging (\ref{17}) into (\ref{15}), one gets
\begin{equation*}
\Gamma _1  \ge \frac{{8(m - p)(\Delta  + \delta _1 )\sqrt {\Delta \delta _1 } }}{{\Delta ^2 (\Delta  + \delta _1  + 2\sqrt {\Delta \delta _1 } )^2 }}(F + 2m\Delta ^2  - p(3\Delta ^2  + 1)),
\end{equation*}
which implies that
\begin{equation*}
AG(G) = \frac{1}{2}\sum\limits_{i \sim j,d_j  = 1} {\left( {\sqrt {\frac{{d_i }}{{d_j }}}  + \sqrt {\frac{{d_j }}{{d_i }}} } \right)}  + \frac{1}{2}\sum\limits_{i \sim j,d_i ,d_j  > 1} {\left( {\sqrt {\frac{{d_i }}{{d_j }}}  + \sqrt {\frac{{d_j }}{{d_i }}} } \right)}~~~~~~~~~~~~~~~~~~~~~~~~~~~~~~~~~~~
\end{equation*}
\begin{equation}\label{18}
~~~~~~ \ge \frac{{p(\delta _1  + 1)}}{{2\sqrt {\delta _1 } }} + \frac{{\sqrt {2(m - p)(\Delta  + \delta _1 )\sqrt {\Delta \delta _1 } } }}{{\Delta (\Delta  + \delta _1  + 2\sqrt {\Delta \delta _1 } )}}\sqrt {F + 2m\Delta ^2  - p(3\Delta ^2  + 1)}.
\end{equation}

Now assume that the equality holds in (\ref{14}). Then all inequalities in above proof
must be equalities. From the equality in (\ref{16}), we have $d_i=d_j=\Delta$ for any non-pendent edge $i\sim j$. It follows from the equality in (\ref{18}) that ${d_i} =
\delta_1$ and $d_j=1$ for any pendent edge $i\sim j$ with pendent vertex $j$.

Next, one has to keep in mind that $G$ is connected. If $p=0$, that is, $G$ has no pendent edge, then $G$ is isomorphic to a $\Delta$-regular graph. If $m=p$, that is, each one of edges in $G$ is pendent edge, then $G$ is isomorphic to ${K_{1,n - 1}}$. Otherwise, $0<p<m$, which implies that $G$ is isomorphic to a $(\Delta,1)$-semiregular graph.

Conversely, it is easy to see that the equality holds in (\ref{14})
for $K_{1,n-1}$ or a regular graph or a $(\Delta,1)$-semiregular graph.
\end{proof}

\noindent{\bfseries Corollary 4.} Let $G$ be a connected graph of
order $n$ with size $m$, maximum degree $\Delta $, minimum degree $\delta$. If $G$ has no pendent vertices, then
\begin{equation*}
AG(G) \ge \frac{{\sqrt {2m(\Delta  + \delta _1 )\sqrt {\Delta \delta _1 } } }}{{\Delta (\Delta  + \delta _1  + 2\sqrt {\Delta \delta _1 } )}}\sqrt {F + 2m\Delta ^2 },
\end{equation*}
with equality if and only if $G$ is isomorphic to a regular graph.\\

Similar to the proof of Theorem 4, we may obtain the following theorem.\\

\noindent{\bfseries Theorem 5.} If $G$ is a connected graph of
order $n$ with size $m$, maximum degree $\Delta $, minimum non-pendent vertex degree $\delta_1$, then
\begin{equation*}
AG(G)  \ge \frac{{p(\delta _1  + 1)}}{{2\sqrt {\delta _1 } }} + \frac{{\sqrt {2(m - p)(\Delta  + \delta _1 )\sqrt {\Delta \delta _1 } } }}{{ \Delta  + \delta _1  + 2\sqrt {\Delta \delta _1 } }}\sqrt {SDD- p(\Delta +\frac{ 1}{\Delta})+2(m-p)},
\end{equation*}
with equality if and only if $G$ is isomorphic to ${K_{1,n - 1}}$ or $G$ is isomorphic to a regular graph or $G$ is isomorphic to a $(\Delta,1)$-semiregular graph.\\

Clearly, $AG(G) \ge\frac{M_1}{2\Delta}$. Here we shall give a minor improvement on this lower bound as follow.\\

\noindent{\bfseries Theorem 6.} If $G$ is a connected graph of
order $n$ with size $m$, maximum degree $\Delta $, minimum non-pendent vertex degree $\delta_1$, then
\begin{equation*}
AG(G) \ge \frac{{p(\delta _1  + 1)}}{{2\sqrt {\delta _1 } }} + \frac{1}{{2\Delta }}(M_1  - p(\Delta  + 1)),
\end{equation*}
with equality if and only if $G$ is isomorphic to ${K_{1,n - 1}}$ or $G$ is isomorphic to a regular graph or $G$ is isomorphic to a $(\Delta,1)$-semiregular graph.

\begin{proof}
It is easy to verify that
\[
\begin{array}{l}
 AG(G) = \frac{1}{2}\sum\limits_{i \sim j,d_j  = 1} {\left( {\sqrt {\frac{{d_i }}{{d_j }}}  + \sqrt {\frac{{d_j }}{{d_i }}} } \right)}  + \frac{1}{2}\sum\limits_{i \sim j,d_i ,d_j  > 1} {\left( {\sqrt {\frac{{d_i }}{{d_j }}}  + \sqrt {\frac{{d_j }}{{d_i }}} } \right)}  \\
  ~~~~~~~~~~\ge \frac{p}{2}\left( {\sqrt {\delta _1 }  + \frac{1}{{\sqrt {\delta _1 } }}} \right) + \frac{1}{{2\Delta }}\left( {\sum\limits_{i \sim j} {(d_i  + d_j )}  - \sum\limits_{i \sim j,d_j  = 1} {(d_i  + d_j )} } \right) \\
 ~~~~~~~~~~ \ge \frac{{p(\delta _1  + 1)}}{{2\sqrt {\delta _1 } }} + \frac{1}{{2\Delta }}(M_1  - p(\Delta  + 1)), \\
 \end{array}
\]
with equality if and only if $G$ has same degree for all non-pendent vertex. The rest of the proof is similar to that of Theorem 4, omitted.
\end{proof}

\section*{3. Comparison between arithmetic-geometric index and other topological indices}

\noindent {\bfseries Theorem 7.} Let $G$ be a connected graph of order $n$, with minimum degree $\delta$.
Then
\begin{equation}\label{19}
GA(G)\leq AG(G)\leq\frac{(\delta+n-1)^2}{4\delta(n-1)}GA(G)
\end{equation}
with left-hand side of equality if and only if $G$ is a regular graph, and right-hand side of equality if
and only if $G$ is isomorphic to ${K_{1,n - 1}}$ or $G$ is isomorphic to $K_n$.

\begin{proof}
Consider the following function
\begin{equation*}
f(x,y) = \frac{{\frac{1}{2}\left( {\sqrt {\frac{x}{y}}  + \sqrt {\frac{y}{x}} } \right)}}{{\frac{{2\sqrt {xy} }}{{x + y}}}} = \frac{{(x + y)^2 }}{{4xy}},
\end{equation*}
where $1\leq\delta\leq x\leq y\leq n-1$. Now, by a simple computation, we get
\[
\frac{{\partial f}}{{\partial x}} = \frac{{4y(x^2  - y^2 )}}{{16x^2 y^2 }} \le 0,
\]
which implies that $f(x,y)$ is decreasing in $x$. Thus, $f(x,y)$ attains the maximum at $(\delta, y)$
for some $\delta\leq y\leq n-1$. On the other hand, it is easy to verify that $f(\delta,y)$ is an increasing function for $y\geq \delta\geq1$.
Therefore,
\[f(x,y)\leq f(\delta,n-1)=\frac{(\delta+n-1)^2}{4\delta(n-1)},\]
which implies that
\begin{equation*}
AG(G)\leq\frac{(\delta+n-1)^2}{4\delta(n-1)}GA(G)
\end{equation*}
with equality if and only if $(d_i,d_j)=(\delta, n-1)$ for every edge $i\sim j$ of $G$.
If $\delta=1$, then $G$ is isomorphic to ${K_{1,n - 1}}$. Otherwise, $\delta\geq 2$,
this time $G$ has no pendent edge. Without loss of generality, suppose that $d_i=\delta$,
then the vertex $i$ has at least two adjacent vertices with degree $n-1$. This implies that $\delta=n-1$.
Therefore, $G$ is isomorphic to ${K_n}$.

The left-hand side of inequality in (\ref{19}) is clearly true. The proof is complete.
\end{proof}

\noindent {\bfseries Corollary 5.} Let $G$ be a connected graph of order $n\geq2$.
Then
\begin{equation*}
AG(G)\leq\frac{n^2}{4(n-1)}GA(G)
\end{equation*}
with equality if and only if $G$ is isomorphic to ${K_{1,n - 1}}$.\\

Denote the chromatic number of a graph $G$ by $\chi(G)$. It was proved in [1] that if $G$ is a connected graph with $\delta\geq2$, then
$\chi(G)\leq\frac{2}{\delta}GA(G)$ with equality if and only if $G$ is isomorphic to ${K_n}$. In [2],
Aouchiche and Hansen proposed the following conjecture.\\

\noindent {\bfseries Conjecture 1} [2]. Let $G$ be a connected graph of order $n$ with $m$ edges and average degree $\overline{d}$.
Then
\begin{equation*}
\chi(G)\leq\frac{2 GA(G)}{\overline{d}}
\end{equation*}
with equality if and only if $G$ is isomorphic to $K_n$.\\

It is easy to see that Conjecture 1 holds for a regular graph $G$ or complete bipartite graph $K_{n_1,n_2}$ of order $n=n_1+n_2$.
Denote the join of $G_1$ and $G_2$ by $G_1\vee G_2$, we define $L(n,k)=K_k\vee \overline{K_{n-k}}$, where $\overline{K_{n-k}}$ is the complement of the complete graph $K_{n-k}$.
Notice that $L(n,1)=K_{1,n-1}$ and $L(n,n-1)=K_n$. Next we assume that $2\leq k\leq n-2$. Clearly, $\chi (L(n,k))=k+1$. By a simple computation, we obtain
\[
\frac{{2GA(L(n,k))}}{{\overline d }} = n \cdot \frac{{\left( {\begin{array}{*{20}c}
   k  \\
   2  \\
\end{array}} \right) + k(n - k)\frac{{2\sqrt {k(n - 1)} }}{{n + k - 1}}}}{{\left( {\begin{array}{*{20}c}
   k  \\
   2  \\
\end{array}} \right) + k(n - k)}} = O(\sqrt n ).
\]
Hence, we arrive at\\

\noindent {\bfseries Theorem 8.} For a fixed number $k$ and sufficiently large $n$, we have
\[
\chi (L(n,k)) \le \frac{{2GA(L(n,k))}}{{\overline d }}.
\]

As we all know, Conjecture 1 is still open. However, if the first geometric-arithmetic index $GA(G)$ is
replaced by arithmetic-geometric index $AG(G)$ in above conjecture, then
 \[
 \chi(G)\leq \frac{2m}{\overline{d}}\leq \frac{2AG(G)}{\overline{d}},
\]
with equality if and only if $G$ is isomorphic to $K_n$.

In [2], it is also pointed out that, there exist graphs with $\chi(G)>\frac{2GA(G)}{\Delta}$. But, similar to Theorem 8,
we easily prove that, for a sufficiently large $n$,
\[
\chi (L(n,k)) \le \frac{{2AG(L(n,k))}}{n-1}.
\]
Thus, the following problem arises: does there exist a graph $G$ satisfying $\chi(G)>\frac{2AG(G)}{\Delta}$?

In the following, we shall consider relations between arithmetic-geometric index $AG(G)$ and atom-bond connectivity index $ABC(G)$.
Let $T^*$ denote the tree obtained by joining the central vertices of two copies of $K_{1,3}$ by an edge.
Das and Trinajsti\'{c} [8] proved that if $G$ is a connected graph with $\Delta-\delta\leq 3$ and it is neither isomorphic to
$K_{1,4}$ nor $T^*$, then $GA(G)>ABC(G)$. Note that $AG(K_{1,4})>ABC(K_{1,4})$ and $AG(T^*)>ABC(T^*)$.
Thus, combining these results with Theorem 7, one gets $AG(G)>ABC(G)$
for any connected graph $G$ with $\Delta-\delta\leq 3$. Next, we give an improvement on this result.\\

\noindent {\bfseries Theorem 9.} Let $G$ be a connected graph of order $n$, with minimum degree $\delta\geq 2$.
Then
\begin{equation}\label{20}
\frac{\delta }{{\sqrt {2\delta  - 2} }}ABC(G) \le AG(G) \le \frac{{n - 1}}{{\sqrt {2n - 4} }}ABC(G).
\end{equation}
Moreover, the left-hand side of equality holds in (\ref{20}) if and only if $G$ is a $\delta$-regular graph, and right-hand side of equality holds in (\ref{20}) if
and only if $G$ is isomorphic to $K_n$.

\begin{proof}
Consider the following function
\begin{equation*}
f(x,y) =\left( \frac{{\frac{1}{2}\left( {\sqrt {\frac{x}{y}}  + \sqrt {\frac{y}{x}} } \right)}}{{\frac{{\sqrt {x+y-2} }}{{\sqrt{x y}}}}}\right)^2 = \frac{{(x + y)^2 }}{{4(x+y-2)}},
\end{equation*}
where $2\leq\delta\leq x\leq y\leq n-1$. Now, by a simple computation, we get
\[
\frac{{\partial f}}{{\partial x}} = \frac{{(x + y)(x + y - 4)}}{{4(x + y - 2)^2 }} \ge 0,
\]
which implies that $f(x,y)$ is increasing in $x$. Thus, $f(x,y)$ attains the minimum at $(\delta, y_1)$
for some $\delta\leq y_1\leq n-1$ and maximum at $(y_2, y_2)$
for some $\delta\leq y_2\leq n-1$ On the other hand, it is easy to verify that $f(\delta,y)$ is an increasing function for $y\geq \delta\geq2$.
Thus,
\[
f(\delta ,\delta ) \le f(x,y) \le f(n - 1,n - 1),
\]
which implies that
\begin{equation*}
\frac{\delta }{{\sqrt {2\delta  - 2} }}ABC(G) \le AG(G) \le \frac{{n - 1}}{{\sqrt {2n - 4} }}ABC(G)
\end{equation*}
with left-hand side of equality if and only if $(d_i,d_j)=(\delta, \delta)$ for every edge $i\sim j$ of $G$,
and right-hand side of equality if and only if $(d_i,d_j)=(n-1,n-1)$ for every edge $i\sim j$ of $G$. Hence, the required result follows.
\end{proof}

Note that it follows from Theorem 9 that, for a graph $G$ with $\delta\geq 2$, $AG(G)>\sqrt{2}ABC(G)$
unless $G$ is isomorphic to $C_n$.

Using the similar technique to the proof in Theorem 9, we easily obtain the following bounds
for the arithmetic-geometric index $AG(G)$ in terms of the symmetric division deg index $SDD(G)$ (the details is omitted).\\

\noindent {\bfseries Theorem 10.} Let $G$ be a connected graph of order $n$, with minimum degree $\delta$.
Then
\begin{equation}\label{21}
\frac{{(\delta  + n - 1)\sqrt {\delta (n - 1)} }}{{2(\delta ^2  + (n - 1)^2) }}SDD(G) \le AG(G) \le \frac{1}{2}SDD(G).
\end{equation}
Moreover, the left-hand side of equality holds in (\ref{21}) if and only if $G$ is isomorphic to ${K_{1,n - 1}}$ or $G$ is isomorphic to $K_n$,
and right-hand side of equality holds in (\ref{21}) if and only if $G$ is a $\delta$-regular graph.

\section*{4. Effect on $GA$ and $AG$ indices of deleting an edge from a graph}

\indent \indent In this section, we mainly discuss the effect on $GA$
and $AG$ indices when an edge is deleted from a graph $G$. First note that $GA$ and $AG$ indices will always decrease when
an edge $e=v_iv_j$, with $d_i=d_j=1$, is deleted from $G$. For simplicity sake, assume that $e=v_iv_j$ is an edge
with non-pendent vertex $v_j$ throughout this section.

\subsection*{4.1. Effect on $GA$ index of deleting an edge}

\indent \indent In [6], Das et al. presented a sufficient condition with $GA(G+e)>GA(G)$
when a new edge $e$ is inserted into the graph $G$. Recently, Chen and Wu [3] pointed out that the result obtained in [6] is not complete.
Furthermore, they established Bollob\'{a}s-Erd\H{o}s-type theorem for the first geometric-arithmetic index of a graph $G$ as follows.\\

\noindent {\bfseries Theorem 11} [3]. Let $G$ be a simple graph with an edge $e=v_iv_j$. Also let $d_r=\max\{d_k|v_iv_k\in E(G)\}$
and $d_s=\max\{d_l|v_jv_l\in E(G)\}$. If one of the following conditions is satisfied, then $GA(G)>GA(G-e)$:
\begin{enumerate}
    \item[(i)] $\max\{\frac{d_i}{d_r},\frac{d_j}{d_s}\}\leq 1,$ or
    \item[(ii)] $\max\{\frac{d_i}{d_j},\frac{d_j}{d_i}\}\leq \min\{\frac{d_i}{d_r},\frac{d_j}{d_s}\}$.
\end{enumerate}

\noindent {\bfseries Example 1.} Let $G$ be the graph as showed in Figure 1, where $d_i=10$ and $d_j=d_r=d_s=1000$. Clearly, $G$ satisfies the condition
(i) of Theorem 11, that is, $\max\{\frac{10}{1000},\frac{1000}{1000}\}\leq 1$, but $GA(G)-GA(G-e)=-0.0447$.

\begin{center}
\includegraphics[width=6cm]{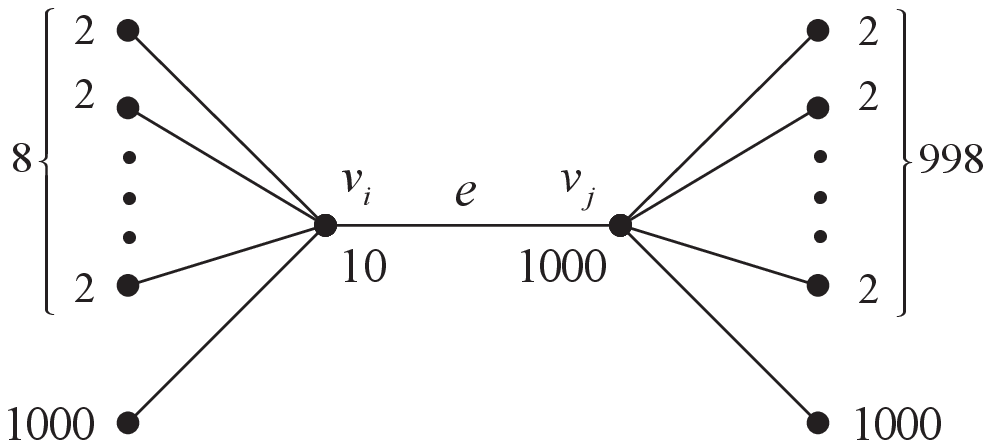}\bigskip
\center{Figure 1: A counterexample to the (i) of Theorem 11.}
\end{center}

\noindent {\bfseries Example 2.} Let $G$ be the graph as showed in Figure 2, where $d_i=100$, $d_j=500$, $d_r=500$ and $d_s=100$.
By a simple calculation, one can see that $GA(G)-GA(G-e)=0.5501$, despite $\max\{\frac{100}{500},\frac{500}{100}\}> \min\{\frac{100}{500},\frac{500}{100}\}$.
Therefore, the (ii) of Theorem 11 is invalid for this graph $G$.

\begin{center}
\includegraphics[width=8cm]{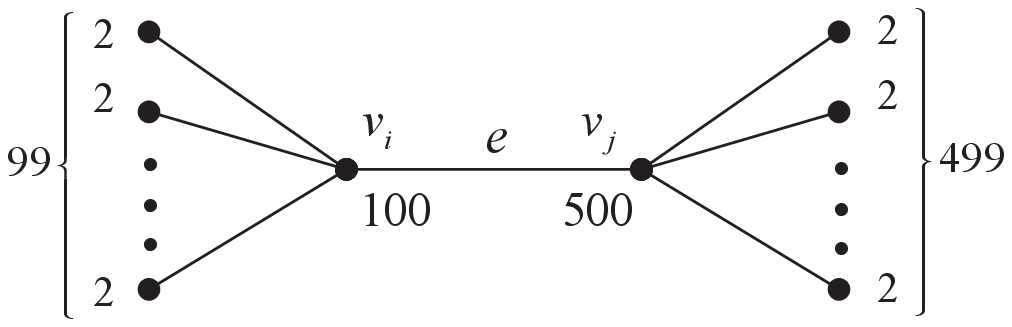}\bigskip
\center{Figure 2: The graph $G$ in Example 2.}
\end{center}

\noindent {\bfseries Example 3.} For given two graphs $G_1$, $G_2$ in Figure 3, one can see that $GA(G_1)-GA(G_1-e)=0.0652$, whereas $GA(G_2)-GA(G_2-e)=-0.0363$.
This example shows that $GA$ index may either increase or decrease when a pendent edge $e$ is deleted from a graph.

\begin{center}
\includegraphics[width=12cm]{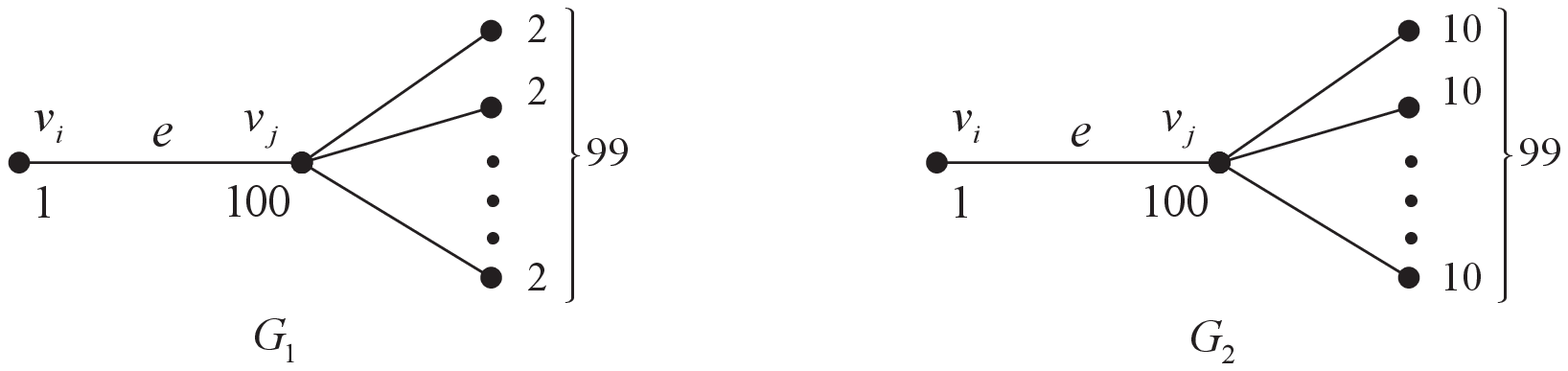}\bigskip
\end{center}
{Figure 3: In the case of $G_1$ the $GA$ index decreases, whereas in the case of $G_2$ the $GA$ index increases.}\\

Assume that $G$ is a simple graph and $e=v_iv_j$ is an edge of $G$ with non-pendent vertex $v_j$. For convenience, define
\[
d_{\min }^{(j)}  = \min \{ d_k |v_k  \in N(v_j )\backslash \{ v_i \} \}~~\text{and}~~ d_{\max }^{(j)}  = \max \{ d_k |v_k  \in N(v_j )\backslash \{ v_i \} \}.
\]
Note that one may give similar definitions when $v_i$ also is a non-pendent vertex of $G$.
Next, by many helpful techniques provided in [3,6] and some analysis, we first provide a sufficient condition for $GA(G)>GA(G-e)$ when $e=v_iv_j$ is a pendent edge of $G$.\\

\noindent {\bfseries Theorem 12.} Assume that $G$ is a simple graph and $e=v_iv_j$ is a pendent edge of $G$ with non-pendent vertex $v_j$.
 If one of the following conditions is satisfied, then $GA(G)>GA(G-e)$:
\begin{enumerate}
    \item[(i)] $d_{\min }^{(j)} \geq d_j$, or
    \item[(ii)] $\frac{\sqrt{d_{\max }^{(j)} }}{2\sqrt{d_j-\frac{1}{2}+6d_{\max }^{(j)} }}\leq \frac{\sqrt{d_j}}{d_j+1}$.
\end{enumerate}

\begin{proof}
Since $e=v_iv_j$ is a pendent edge of $G$ with non-pendent vertex $v_j$, then
\begin{equation}\label{22}
GA(G) - GA(G - e) = 2\sum\limits_{v_k  \in N(v_j )\backslash \{ v_i \} } {\left( {\frac{{\sqrt {d_j d_k } }}{{d_j  + d_k }} - \frac{{\sqrt {(d_j  - 1)d_k } }}{{d_j  + d_k  - 1}}} \right)}  + 2\frac{{\sqrt {d_j } }}{{d_j  + 1}}.
\end{equation}

If $d_{\min }^{(j)} \geq d_j$, then $d_k\geq d_j$ for any $v_k  \in N(v_j )\backslash \{ v_i \}$. Notice that $f(x)=\frac{1}{x+\frac{1}{x}}$ is an increasing function for $x\in(0,1]$.
Thus one easily see that
\begin{equation}\label{23}
\frac{{\sqrt {(d_j  - 1)d_k } }}{{d_j  + d_k  - 1}} - \frac{{\sqrt {d_j d_k } }}{{d_j  + d_k }} < 0,
\end{equation}
which implies that $GA(G)>GA(G-e)$. Hence, the (i) follows.

Now suppose that $d_k\leq d_j-1$ for some $v_k  \in N(v_j )\backslash \{ v_i \}$. Then
\begin{equation*}
\frac{{\sqrt {(d_j  - 1)d_k } }}{{d_j  + d_k  - 1}} - \frac{{\sqrt {d_j d_k } }}{{d_j  + d_k }}~~~~~~~~~~~~~~~~~~~~~~~~~~~~~~~~~~~~~~~~~~~~~~~~~~~~~~~~~~~~~~~~~~~~~~~~~~~~~~~~~~~~~~~~~~~~~~~~~~~~~~~~~~~~~~~~~~~~~~~~~~~~~~
\end{equation*}
\begin{equation*}
 = \frac{{\sqrt {d_k } }}{{(d_j  + d_k )(d_j  + d_k  - 1)}}\left( {\sqrt {d_j  - 1} (d_j  + d_k ) - \sqrt {d_j } (d_j  + d_k  - 1)} \right)
\end{equation*}
\begin{equation*}\small
 \le \frac{{\sqrt {d_k } \left( {\sqrt {d_j  - 1} (d_j  + d_k ) - (d_j  + d_k  - 1)(\sqrt {d_j  - 1}  + \frac{1}{{2\sqrt {d_j  - 1} }} - \frac{1}{{8(d_j  - 1)^{3/2} }})} \right)}}{{(d_j  + d_k )(d_j  + d_k  - 1)}}
\end{equation*}
\begin{equation*}
 = \frac{{\sqrt {d_k } }}{{(d_j  + d_k )(d_j  + d_k  - 1)}}\left( {d_j  - d_k  - 1 + \frac{{d_j  + d_k  - 1}}{{4(d_j  - 1)}}} \right)\frac{1}{{2\sqrt {d_j  - 1} }}
\end{equation*}
\begin{equation*}
 \le \frac{{\sqrt {d_k } }}{{2(d_j  - 1)\sqrt {d_j  - \frac{1}{2} + 6d_k } }} \cdot \frac{{(d_j  - d_k  - \frac{1}{2})\sqrt {d_j  - 1} \sqrt {d_j  - \frac{1}{2} + 6d_k } }}{{(d_j  + d_k )(d_j  + d_k  - 1)}}
\end{equation*}
\begin{equation*}
 < \frac{{\sqrt {d_k } }}{{2(d_j  - 1)\sqrt {d_j  - \frac{1}{2} + 6d_k } }} \cdot \frac{{(d_j  - d_k  - \frac{1}{2})(d_j  + 3d_k  - \frac{1}{2})}}{{(d_j  + d_k )(d_j  + d_k  - 1)}}
\end{equation*}
\begin{equation*}
 < \frac{{\sqrt {d_k } }}{{2(d_j  - 1)\sqrt {d_j  - \frac{1}{2} + 6d_k } }}
\end{equation*}
\begin{equation}\label{24}
\leq \frac{{\sqrt {d_{\max }^{(j)} } }}{{2(d_j  - 1)\sqrt {d_j  - \frac{1}{2} + 6d_{\max }^{(j)} } }},
\end{equation}
where the last inequality holds as $d_{\max }^{(j)}\geq d_k$ and $f(x)=\frac{\sqrt{x}}{\sqrt{a+6x}}$ is an increasing function in $x>0$ whenever $a>0$. Now, using (\ref{23}) and (\ref{24}), one may get
\begin{equation}\label{25}
\sum\limits_{v_k  \in N(v_j )\backslash \{ v_i \} } {\left( {\frac{{\sqrt {(d_j  - 1)d_k } }}{{d_j  + d_k  - 1}} - \frac{{\sqrt {d_j d_k } }}{{d_j  + d_k }}} \right)}  < \frac{{\sqrt {d_{\max }^{(j)} } }}{{2\sqrt {d_j  - \frac{1}{2} + 6d_{\max }^{(j)} } }}.
\end{equation}
Therefore, if the condition (ii) is satisfied, then $GA(G)>GA(G-e)$.

The proof is complete.
\end{proof}

Remark that, in Example 3, it is easy to verify that $G_1$ satisfies the condition (ii) of Theorem 12.
So, $GA(G_1)>GA(G_1-e)$. However, $G_2$ does not satisfy each of conditions of Theorem 11.\\

\noindent {\bfseries Theorem 13.} Assume that $G$ is a graph with non-pendent edge $e=v_iv_j$.
 If one of the following conditions is satisfied, then $GA(G)>GA(G-e)$:
\begin{enumerate}
    \item[(i)] $\max \left\{ {\frac{{d_i }}{{d_{\min }^{(i)} }},\frac{{d_j }}{{d_{\min }^{(j)} }}} \right\} \le 1$, or
    \item[(ii)] $\max \left\{ {\frac{{d_i }}{{d_j }},\frac{{d_j }}{{d_i }}} \right\} \le \min \left\{ {\frac{{d_i  - \frac{1}{2}}}{{d_{\max }^{(i)} }},\frac{{d_j  - \frac{1}{2}}}{{d_{\max }^{(j)} }}} \right\}$.
\end{enumerate}

\begin{proof}
Since $e=v_iv_j$ is a non-pendent edge of $G$. Then, from the definition of $GA$ index, one gets
\begin{equation*}
GA(G) - GA(G - e) = 2\sum\limits_{v_l  \in N(v_i )\backslash \{ v_j \} } {\left( {\frac{{\sqrt {d_i d_l } }}{{d_i  + d_l }} - \frac{{\sqrt {(d_i  - 1)d_l } }}{{d_i  + d_l  - 1}}} \right)}~~~~~~~~~~~~~~~~~~~~~~~~~~~~~~~~~~~
\end{equation*}
\begin{equation}\label{26}
~~~~~~~~~~~~~~~~~~~~~~~~~~~~~~ + 2\sum\limits_{v_k  \in N(v_j )\backslash \{ v_i \} } {\left( {\frac{{\sqrt {d_j d_k } }}{{d_j  + d_k }} - \frac{{\sqrt {(d_j  - 1)d_k } }}{{d_j  + d_k - 1}}} \right)}  + 2\frac{{\sqrt {d_i d_j } }}{{d_i  + d_j }}.
\end{equation}

If $G$ satisfies the condition (i), then $d_l\geq d_i$ for any $v_l  \in N(v_i )\backslash \{ v_j \}$ and $d_k\geq d_j$ for any $v_k  \in N(v_j )\backslash \{ v_i \}$.
Now, from (\ref{23}), one can easily see that
\begin{equation*}
\frac{{\sqrt {(d_i  - 1)d_l } }}{{d_i  + d_l  - 1}} - \frac{{\sqrt {d_i d_l } }}{{d_i  + d_l }} < 0,
\end{equation*}
and
\begin{equation*}
\frac{{\sqrt {(d_j  - 1)d_k } }}{{d_j  + d_k  - 1}} - \frac{{\sqrt {d_j d_k } }}{{d_j  + d_k }} < 0.
\end{equation*}
So, it follows from (\ref{26}) that $GA(G)>GA(G-e)$. Otherwise, again from (\ref{23}) and (\ref{24}), one has
\begin{equation*}
\sum\limits_{v_k  \in N(v_j )\backslash \{ v_i \} } {\left( {\frac{{\sqrt {(d_j  - 1)d_k } }}{{d_j  + d_k  - 1}} - \frac{{\sqrt {d_j d_k } }}{{d_j  + d_k }}} \right)}  < \frac{{\sqrt {d_{\max }^{(j)} } }}{{2\sqrt {d_j  - \frac{1}{2} + 6d_{\max }^{(j)} } }}.
\end{equation*}
Similarly,
\begin{equation*}
\sum\limits_{v_l  \in N(v_i )\backslash \{ v_j \} } {\left( {\frac{{\sqrt {(d_i  - 1)d_l } }}{{d_i  + d_l  - 1}} - \frac{{\sqrt {d_i d_l } }}{{d_i  + d_l }}} \right)}  < \frac{{\sqrt {d_{\max }^{(i)} } }}{{2\sqrt {d_i  - \frac{1}{2} + 6d_{\max }^{(i)} } }}.
\end{equation*}

The following proof is similar to that of Theorem 3.4 in [3]. For the convenience of readers, here we give the detailed proof. Let $t_1=\max \left\{ {\frac{{d_i }}{{d_j }},\frac{{d_j }}{{d_i }}} \right\}$, $t_2=\frac{{d_i  - \frac{1}{2}}}{{d_{\max }^{(i)} }}$ and $t_3=\frac{{d_j  - \frac{1}{2}}}{{d_{\max }^{(j)} }}$. Without loss of generality, assume that $t_2\leq t_3$. The condition (ii) implies that $1\leq t_1\leq t_2\leq t_3$.
After some rearrangements, one has
\[
\frac{{\sqrt {d_{\max }^{(i)} } }}{{2\sqrt {d_i  - \frac{1}{2} + 6d_{\max }^{(i)} } }} + \frac{{\sqrt {d_{\max }^{(j)} } }}{{2\sqrt {d_j  - \frac{1}{2} + 6d_{\max }^{(j)} } }} = \frac{1}{{2\sqrt {t_2  + 6} }} + \frac{1}{{2\sqrt {t_3  + 6} }}~~~~~~~~~~~~~~~~~~~~~~~~~~~~~~~~~~~~~~~~~~~~~~~~~~~~~~~~~~
\]
\[
~~~~~~~~~~~~~~~~~~~~~~~~~~~~~~~~~~~~~~~~~~~ \le \frac{1}{{\sqrt {t_2  + 6} }} \le \frac{{\sqrt {t_1 } }}{{t_1  + 1}} = \frac{{\sqrt {d_i d_j } }}{{d_i  + d_j }}.
\]
Hence, it follows from (\ref{26}) that $GA(G)>GA(G-e)$.

The proof is complete.
\end{proof}

Remark that, in Example 2, $d_i=100$, $d_j=500$, $d_{\max }^{(i)}=d_{\max }^{(j)}=2$. Clearly, $G_1$ satisfies the condition (ii) of Theorem 13.
So, $GA(G)>GA(G-e)$. Hence, Theorem 13 is an improvement on Theorem 11. In addition, if $G$ has an edge $e=v_iv_j$ with the property (i) in Theorem 13, we say $e$ is an \emph{ascending edge} of $G$.\\

\noindent {\bfseries Corollary 6.} If $e=v_iv_j$ is an ascending edge of $G$, then $GA(G)>GA(G-e)$.

\subsection*{4.2. Effect on $AG$ index of deleting an edge}

\noindent {\bfseries Theorem 14.} Let $e=v_iv_j$ be an edge of a graph $G$ with non-pendent vertex $v_j$.
 If one of the following conditions is satisfied, then $AG(G)>AG(G-e)$:
\begin{enumerate}
    \item[(i)] $\min \left\{ {\frac{{d_i }}{{d_{\max }^{(i)} }},\frac{{d_j }}{{d_{\max }^{(j)} }}} \right\} > 1$, or
    \item[(ii)] $\frac{{d_{\max }^{(i)}  - d_i  + 1}}{{2\sqrt {d_{\max }^{(i)} } \sqrt {d_i } }} + \frac{{d_{\max }^{(j)}  - d_j  + 1}}{{2\sqrt {d_{\max }^{(j)} } \sqrt {d_j } }} \le \frac{{d_i  + d_j }}{{\sqrt {d_i d_j } }},$
\end{enumerate}
where $d_i/d_{\max }^{(i)}$ is stipulated as $\propto$ and $\frac{{d_{\max }^{(i)}  - d_i  + 1}}{{2\sqrt {d_{\max }^{(i)} } \sqrt {d_i } }} =\frac{d_{\max }^{(i)} }{2}=0$ when $v_i$ is a pendent vertex.

\begin{proof}
First suppose that $e=v_iv_j$ is non-pendent edge of $G$. Then, in the light of the definition of $GA$ index,
\[
AG(G) - AG(G - e) = \frac{1}{2}\sum\limits_{v_k  \in N(v_i )\backslash \{ v_j \} } {\left( {\frac{{d_i  + d_k }}{{\sqrt {d_i d_k } }} - \frac{{d_i  + d_k  - 1}}{{\sqrt {(d_i  - 1)d_k } }}} \right)}~~~~~~~~~~~~~~~~~~~~~~~~~~~~~~~
\]
\[
 ~~~~~~~~~~~~~~~~~~~~~+ \frac{1}{2}\sum\limits_{v_l  \in N(v_j )\backslash \{ v_i \} } {\left( {\frac{{d_j  + d_l }}{{\sqrt {d_j d_l } }} - \frac{{d_j  + d_l  - 1}}{{\sqrt {(d_j  - 1)d_l } }}} \right)}  + \frac{{d_i  + d_j }}{{2\sqrt {d_i d_j } }}.
\]

If $G$ satisfies the condition (i), then $d_i> d_k$ for any $v_k  \in N(v_i )\backslash \{ v_j \}$ and $d_j>d_l$ for any $v_l  \in N(v_j )\backslash \{ v_i \}$.
Since $f(x)=x+\frac{1}{x}$ is an increasing function for $x\geq 1$. Thus,
\begin{equation}\label{27}
\frac{{d_i  + d_k }}{{\sqrt {d_i d_k } }} > \frac{{d_i  + d_k  - 1}}{{\sqrt {(d_i  - 1)d_k } }}
\end{equation}
and
\[
\frac{{d_j  + d_l }}{{\sqrt {d_j d_l } }} > \frac{{d_j  + d_l  - 1}}{{\sqrt {(d_j  - 1)d_l } }}.
\]
Hence, $AG(G)>AG(G-e)$.

If $d_k\geq d_i>1$ for some $v_k  \in N(v_i )\backslash \{ v_j \}$, then
\[
\frac{{d_i  + d_k  - 1}}{{\sqrt {(d_i  - 1)d_k } }} - \frac{{d_i  + d_k }}{{\sqrt {d_i d_k } }} = \frac{{\sqrt {d_i } \sqrt {(d_i  - 1)} (d_i  + d_k  - 1) - (d_i  - 1)(d_i  + d_k )}}{{(d_i  - 1)\sqrt {d_k } \sqrt {d_i } }}~~~~~~~~~~~~~~~
\]
\[
~~~~~~~~~~~~~~~~~~~~~~~< \frac{{(d_i  - \frac{1}{2})(d_i  + d_k  - 1) - (d_i  - 1)(d_i  + d_k )}}{{(d_i  - 1)\sqrt {d_k } \sqrt {d_i } }}
\]
\[
~~~~~~~~~~~~~~~~~~~~~~~ = \frac{{d_k  - d_i  + 1}}{{2(d_i  - 1)\sqrt {d_k } \sqrt {d_i } }}
\]
\begin{equation}\label{28}
~~~~~~~~~~~~~~~~~~~~~~~ \le \frac{{d_{\max }^{(i)}  - d_i  + 1}}{{2(d_i  - 1)\sqrt {d_{\max }^{(i)} } \sqrt {d_i } }},
\end{equation}
where the last inequality holds as $d_{\max }^{(i)}\geq d_k$ and $f(x)=\frac{x-a}{b\sqrt{x}}$ $(a,b>0)$ is an increasing function for $x\geq 0$. Hence, from (\ref{27}) and (\ref{28}), one has
\[
\sum\limits_{v_k  \in N(v_i )\backslash \{ v_j \} } {\left( {\frac{{d_i  + d_k  - 1}}{{\sqrt {(d_i  - 1)d_k } }} - \frac{{d_i  + d_k }}{{\sqrt {d_i d_k } }}} \right)}  < \frac{{d_{\max }^{(i)}  - d_i  + 1}}{{2\sqrt {d_{\max }^{(i)} } \sqrt {d_i } }}.
\]
Similarly,
\[
\sum\limits_{v_l  \in N(v_j )\backslash \{ v_i \} } {\left( \frac{{d_j  + d_l  - 1}}{{\sqrt {(d_j  - 1)d_l } }} -{\frac{{d_j  + d_l }}{{\sqrt {d_j d_l } }}} \right)}  < \frac{{d_{\max }^{(j)}  - d_j  + 1}}{{2\sqrt {d_{\max }^{(j)} } \sqrt {d_j } }}.
\]
Therefore, if $G$ is a graph with the condition (ii), then $AG(G)>AG(G-e)$.

If $e=v_iv_j$ is a pendent edge, then $v_i$ is a pendent vertex. After a simple check, the result still follows.

The proof is complete.
\end{proof}

Let $H$ be any graph of order $n-2$ with maximum degree $\Delta(H)< n-3$ and $G=K_2\vee H$. If $e=v_iv_j$ is the edge with $d_i=d_j=n-1$, then $AG(G)>AG(G-e)$. We say $e=v_iv_j$ is a \emph{descending edge} of $G$ if the edge $e$ has the property (i) in Theorem 14.\\

\noindent {\bfseries Corollary 7.} If $e=v_iv_j$ is a descending edge of $G$, then $AG(G)>AG(G-e)$.\\

\noindent {\bfseries Example 4.} For given two graphs $G_1$, $G_2$ in Figure 4, one can see that $AG(G_1)-AG(G_1-e)=-1.0170$, whereas $AG(G_2)-AG(G_2-e)=0.6309$.
In fact, $G_2$ satisfies the (ii) of Theorem 14.
This example also shows that $AG$ index may either increase or decrease when an ascending edge $e$ is deleted from a graph. However, Corollary 6 implies that $GA$ indices of $G_1$ and $G_2$ are all decrease when the ascending edge $e$ is deleted. So there are considerable differences between $GA$ and $AG$ indices of graphs.

\begin{center}
\includegraphics[width=12cm]{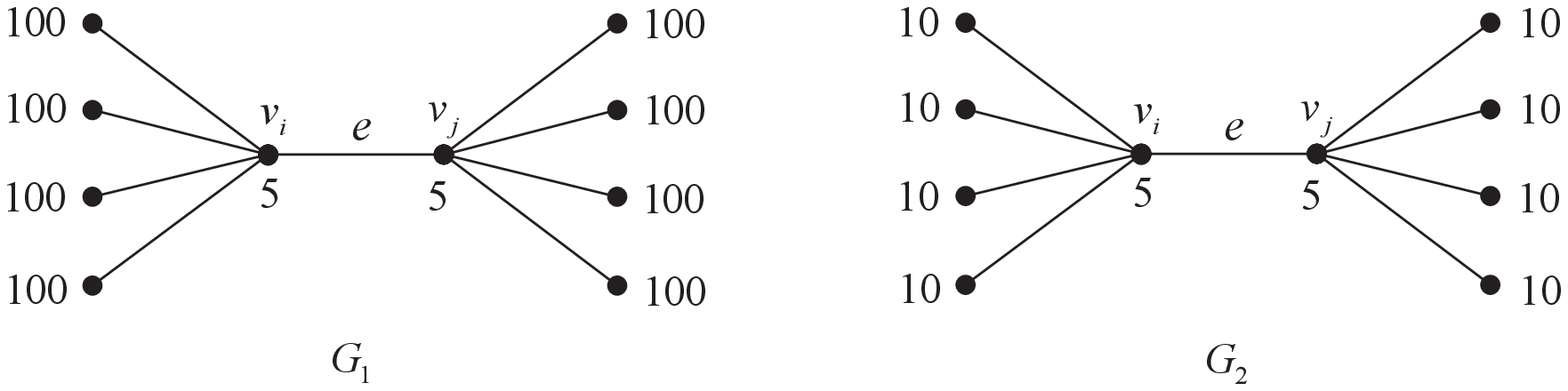}\bigskip
\end{center}
{Figure 4: In the case of $G_1$ the $AG$ index increases, whereas in the case of $G_2$ the $AG$ index decreases.}\\

Finally, we suggest the following problem.\\

\noindent {\bfseries Problem 1:} Is there a graph $G$ such that $GA(G)=GA(G-e)$ or $AG(G)=AG(G-e)$ for some edge $e\in E(G)$?\\

\noindent \textbf{Acknowledgements} This work was in part supported by the National Natural
Science Foundation of China (Nos. 11801521).

\end{document}